\documentclass[11pt]{amsart} 

\usepackage[hyperindex,breaklinks]{hyperref} % useful for arxiv

\usepackage[totalwidth=450pt, totalheight=630pt,centering]{geometry} 
\usepackage{enumitem}

\usepackage[%numbers,
sort&compress]{natbib}

\newenvironment{theorem}[1][Theorem.]{\begin{trivlist}
\item[\hskip \labelsep {\bfseries #1}]\slshape}{\end{trivlist}}

\newtheorem{lemma}{Lemma}
\newenvironment{remarks}[1][Remarks.]{\begin{trivlist}
\item[\hskip \labelsep {\bfseries #1}]}{\end{trivlist}}

\newcommand{\N}{\mathbb{N}}
\newcommand{\R}{\mathbb{R}}
\newcommand{\Z}{\mathbb{Z}} 

\newcommand{\ii}{\mathrm{i}}
\newcommand{\ee}{\mathrm{e}}
\newcommand{\pii}{\mathrm{\pi}}
\newcommand{\dd}{\mathrm{d}}

\newcommand{\Uniform}{\mathrm{U}} 
\newcommand{\uniform}{\mathrm{u}} 

\renewcommand{\phi}{\varphi} 

\newcommand{\1}{\mathbf{1}}           % indicator function

\begin{document}
%%%%%%%%%%%%%%%%%%%%%%%%%%%%%%%%%%%%%%%%%%%%%%%%%%%%%%%%%%%%%%%%%%%%%%
\title[Maximal probabilities of convolutions]{Maximal 
probabilities of convolution powers\\
  of discrete uniform  distributions}

% Information for first author

\author{Lutz Mattner}
\address{ Universit\"at zu L\"ubeck, Institut f\"ur  Mathematik, 
 Wallstr.\ 40, 23560 L\"ubeck, Germany}
\email{mattner@math.uni-luebeck.de}

% Information for second author

\author{Bero Roos} 
\address{Universit\"at Hamburg, Department Mathematik, SPST,
Bundesstr.\ 55, 20146 Hamburg, Germany}
\email{roos@math.uni-hamburg.de}

%%%%%%%%%%%%%%%%%%%%%%%%%%%%%%%%%%%%%%%%%%%%%%%%%%%%%%%%%%%%%%%%%%%%%%
%    General info


\date{June 6, 2007 (filename: \texttt{MaxProbUnifConv20070606D.tex})}
%\date{January 1, 2001 and, in revised form, June 22, 2001.}

%%%%%%%%%%%%%%%%%%%%%%%%%%%%%%%%%%%%%%%%%%%%%%%%%%%%%%%%%%%%%%%%%%%%%%

\begin{abstract} 
We prove optimal constant over root $n$ upper bounds
for  the maximal probabilities of $n$th convolution powers of 
discrete uniform distributions. 
\end{abstract}

\maketitle

%%%%%%%%%%%%%%%%%%%%%%%%%%%%%%%%%%%%%%%%%%%%%%%%%%%%%%%%%%%%%%%%%%%%%%
%\section{Result and remarks}
%%%%%%%%%%%%%%%%%%%%%%%%%%%%%%%%%%%%%%%%%%%%%%%%%%%%%%%%%%%%%%%%%%%%%%
For $\ell, n \in \N := \{1,2,3,\ldots\}$, let  $ \mathrm{U}_\ell^{\ast
n}$ denote the $n$th convolution power of the discrete uniform
distribution $ \mathrm{U}_\ell := \frac
1{\ell}(\delta_0+\ldots+\delta_{\ell-1})$. Let $\mathrm{u}^{\ast
n}_\ell$ denote the density of $\mathrm{U}_\ell^{\ast n}$ with respect
to counting measure. Thus,  writing $\1_A(x) := 1$  if $x\in A$ and
$\1_A(x) := 0$ otherwise, we  have for $\ell \in \N$ and $k\in \Z$
\begin{equation} \label{u1,u2}
 \mathrm{u}^{\ast 1}_\ell(k) \,=\, 
 \frac 1\ell\1_{\{0,\ldots,\ell-1\}}(k),\qquad
 \mathrm{u}^{\ast 2}_\ell(k) 
    \,=\, \frac{\ell -|\ell -1-k| }{\ell^2}
    \1_{\{0,\ldots,2\,(\ell-1)\}}(k)
\end{equation}
and the general formula
\begin{equation*} 
\mathrm{u}^{\ast n}_\ell(k)  = \frac1{\ell^n}
\sum_{j=0}^{\lfloor k/\ell\rfloor}
(-1)^j\binom{n}{j}\binom{n+k-\ell j -1}{n-1}
\quad\qquad (\ell,n\in\N,\, k\in \Z)
\end{equation*}
where $\sum_{j=a}^b := 0$ if $a>b$,   for which we refer to 
\citet[pp.~39--43]{moi1756} or  \citet[pp.~34--35]{hal98}.  
The  purpose of this note is to provide a sharp  upper bound for the
maximal probabilities or concentrations
\begin{equation}\label{Def celln}
  c_{\ell,n} \,:=\,  \max_{k\in \Z}\, \mathrm{u}^{\ast n}_\ell(k)
\end{equation}
of $ \mathrm{U}_\ell^{\ast n}$, see Remarks  \ref{RemApplConIneq}  
and   \ref{RemarkSSHvW} below  for possible  applications. From
\eqref{u1,u2}, we obviously get 
\begin{equation}\label{cln for n=1,2}
 c_{\ell,1} \,=\, c_{\ell,2} \,=\, \frac1{\ell}
 \quad\qquad(\ell \in \N)
\end{equation}
In what follows, we exclude the trivial case of $ \mathrm{U}_1^{\ast
n}= \delta_0$  and hence always assume that $\ell \ge 2$.
%%%%%%%%%%%%%%%%%%%%%%%%%%%%%%%%%%%%%%%%%%%%%%%%%%%%%%%%%%%%%%%%%%%%%%
\begin{theorem}
Let $\ell,n\in \N$ with $\ell \ge 2$ and  let $c_{\ell,n}$ be defined
by \eqref{Def celln}. If $n\neq 2$ or $\ell \in \{2,3,4\}$, then
\begin{equation}\label{MainInequality}
  c_{\ell,n} \,<\, \sqrt{ \frac{6}{\pi(\ell^2-1)n}}
\end{equation}
holds. If $n=2$ and $\ell\ge 5$, then inequality
\eqref{MainInequality} has to be reversed. 
\end{theorem}
%%%%%%%%%%%%%%%%%%%%%%%%%%%%%%%%%%%%%%%%%%%%%%%%%%%%%%%%%%%%%%%%%%%%%%
\begin{remarks} 
\begin{enumerate}[label=\textbf{(\alph*)},ref=(\alph*),fullwidth,
labelsep=1ex,itemindent=\parindent,listparindent=\parindent] 

\item \label{RemOnSharpness} 
Let us fix $\ell \ge 2$ and denote by $\mu := (\ell-1)/2$
and $\sigma^2:=(\ell^2-1)/12$ the mean and the variance of 
$\mathrm{U}_\ell$ and let 
$\phi(x):=(1/\sqrt{2\pii})\exp(-x^2/2)$ for $x\in \R$.  By the
local central limit theorem, see e.g.\ \citet[p.~130]{dur05},  we
then have  $\lim_{n\rightarrow \infty}
\sup_{k\in\Z}|\sqrt{n}\,\mathrm{u}^{\ast n}_\ell(k) - \frac 1\sigma 
\phi\big(  (k-n\mu)/(\sigma\sqrt{n})\big)| = 0$. Since the function
$\phi$ is maximal and continuous at zero, we easily get
$\lim_{n\rightarrow \infty} \sqrt{n}\,c_{\ell,n}=\frac 1\sigma\phi(0)=
\sqrt{6/(\pi(\ell^2-1)) }$. Hence \eqref{MainInequality} is sharp for
$n\rightarrow\infty$ and every $\ell$, in the sense that the quotient 
of both sides of the inequality  converges to one. 

\item A corollary to the theorem is the simpler bound 
\begin{equation}\label{CorToMainIneq}
  c_{\ell,n} \,<\, \frac{2\,\sqrt{2/\pi}}{\ell\,\sqrt{n}}
  \quad\qquad (\ell,n\in \N,\, \ell \ge 2)
\end{equation}
obtained by using   $\ell^2-1 \ge 3 \ell^2/4$ in inequality
\eqref{MainInequality}  if  $n\neq 2$,  and \eqref{cln for n=1,2} for
$n=2$. By the previous Remark \ref{RemOnSharpness} and by comparison
with \eqref{MainInequality}, it is obvious that 
\eqref{CorToMainIneq}  is sharp for $n\rightarrow\infty$ only if 
$\ell=2$.
Inequality \eqref{CorToMainIneq} is contained in \citet[]{bre04}:
His Lemme 33.4.4 a) states, in our notation,
\begin{equation}\label{eq213}
 c_{\ell,n} \,\leq\, \frac 2{\ell}  c_{2,n}   
  \quad\qquad (\ell,n\in \N,\, \ell \ge 2)
\end{equation}
which, by the standard Wallis product inequality recalled in Remark
\ref{Remark on Wallis} below, implies \eqref{CorToMainIneq}. Further,
inequality \eqref{CorToMainIneq} results if  Bretagnolle's
Th\'eor\`eme 33.1.1 is applied to random variables each  with
distribution $ \mathrm{U}_\ell$.

\item The existence of some constant $A<\infty$ with
\begin{equation} \label{IneqConsUnspec}
  c_{\ell,n} \,<\, \frac{A }{\ell\,\sqrt{n}}
  \quad\qquad (\ell,n\in \N,\, \ell \ge 2)
\end{equation} 
already follows from
\citeauthor{kes69}'s \citeyearpar{kes69}
concentration inequality for sums of independent real-valued random
variables and, alternatively, from 
\citeauthor{gam73}'s \citeyearpar{gam73}
sharper result for the special case of identically distributed 
symmetric unimodal lattice random  variables. 
In the case considered here, Gamkrelidze's result yields our 
inequality (\ref{MainInequality})
with an additional $O(n^{-1})$-term on the right-hand side.
For a general introduction to concentration inequalities and further 
results, see \citet[sections 1.5 and  2.4]{pet95}.

\item  \label{RemApplConIneq}
\citet[]{bre04},  \citet[]{rog87}, and  
\citet[in particular Theorem 10 and the unproved 
remark on p.~97]{LR94} state upper
bounds for concentrations of sums of independent real-valued random
variables $X_j$  in terms of concentrations of sums of certain
independent $Y_j$ with distributions $\mathrm{U}_{\ell_j}$. (Both
Bretagnolle and Rogozin refer to an unpublished preprint of
Bretagnolle from 1982. Leader and  Radcliffe fail to give appropriate
references to the probabilistic literature.)  Of these authors only
Bretagnolle goes on to deduce an analytically convenient and still 
rather sharp bound, using in particular inequality  \eqref{eq213}.
Possibly the present asymptotically sharper inequality 
\eqref{MainInequality} could serve to improve  Bretagnolle's 
result.

\item
Since  $\Uniform_\ell^{*n}$ is a convolution of distributions 
unimodal on $\Z$ and with some centers of symmetry, it follows from
the well-known discrete Wintner theorem, see   
\citet[page~109, Theorem~4.7]{DJD88} or, more precisely, 
\citet[Lemma 3.3]{mat06}, that the density $\uniform_\ell^{*n}$ is
maximized at the one or two central points of its support
$\{0,\ldots,n(\ell-1)\}$, so that we have
\begin{equation} \label{WhereUlnMax}
 c_{\ell,n} 
  \,=\, \uniform_\ell^{*n}(\Big\lfloor\frac{n(\ell-1)}{2}\Big\rfloor)
  \,=\, \uniform_\ell^{*n}(\Big\lceil\frac{n(\ell-1)}{2}
  \Big\rceil)
\end{equation}

\item  \label{Remark on Wallis}
For $\ell = 2$, the theorem  reduces to the familiar Wallis 
product inequality for the maximal probabilities of symmetric binomial
distributions,
\begin{equation}\label{WallisIneq}
 \binom{2k}{k}2^{-2k} \,<\, \frac 1{\sqrt{\pii k}} 
 \quad\qquad (k\in\N)
\end{equation}
since 
$c_{2,2k-1}=\binom{2k-1}{k}2^{-(2k-1)}
=\binom{2k}{k}2^{-2k}=c_{2,2k}$,
and since the right-hand side of \eqref{MainInequality} for $\ell =2$
and $n=2k$ or $n=2k-1$ is, respectively,  equal to or greater than the
right-hand side of \eqref{WallisIneq}.

\item A concentration bound related to the present theorem is given
in \citet[]{kan76} and in \citet[]{MR06}.   Theorem 2.1  of
the latter paper specialized to $p_j=2/3$ for every $j$ and the
formulas (15) and (8) there yield the inequalities, 
sharp for $n\rightarrow\infty$,
\begin{equation}  \label{IneqFromMR}
 \max_{k\in\Z}\mathrm{U}_3^{\ast n}(\{k,k+1\}) \,<\, G(2n/3)
 \,<\, \sqrt{\frac{3}{\pii n}}\quad\qquad (n\in\N)
\end{equation} 
where 
$ G(\lambda):= 
\ee^{-\lambda}(\mathrm{I}_0(\lambda)+\mathrm{I}_1(\lambda))$
for $\lambda \in {[0,\infty[}$ and $\mathrm{I}_0, \mathrm{I}_1$ denote
the usual modified Bessel functions. Since the left-hand side of
\eqref{IneqFromMR} is $\leq 2c_{3,n}$, the inequality between the
extreme members of \eqref{IneqFromMR} also follows from  the special
case $\ell =3$ of the present theorem.

\item \label{RemarkSSHvW}
A recent application of upper bounds for $c_{\ell,n}$  occurred in 
the  construction of a two-dimensional transient  but polygonally
recurrent random walk by \citet[]{SW06}, who  proved and used 
\eqref{IneqConsUnspec}, see their Lemmas 6 and 1.
 
\end{enumerate}
\end{remarks}
\medskip
%%%%%%%%%%%%%%%%%%%%%%%%%%%%%%%%%%%%%%%%%%%%%%%%%%%%%%%%%%%%%%%%%%%%%%
We will need two standard lemmas for the proof of the theorem. In what
follows, we use the adjectives ``positive'', ``increasing'' etc.\ in
the wide sense. Thus, e.g., a function $f$ with $0\leq f(x)\leq f(y)$
for $x< y$ is called positive and increasing.
%%%%%%%%%%%%%%%%%%%%%%%%%%%%%%%%%%%%%%%%%%%%%%%%%%%%%%%%%%%%%%%%%%%%%%
\begin{lemma}\label{l3}\slshape
Let $a\in{]0,\infty[}$ and let $f,g\,:\,[-a,a]\rightarrow\R$ be
functions with $f$  even,  $f$ decreasing on $[0,a]$, and  $g$  
convex. Then 
\begin{equation*}
\int_{-a}^a f(x)\,g(x)\,\dd x \,\leq\, \frac{1}{2a}
\int_{-a}^a f(x)\,\dd x \,\int_{-a}^ag(x)\,\dd x
\end{equation*}
\end{lemma}
%%%%%%%%%%%%%%%%%%%%%%%%%%%%%%%%%%%%%%%%%%%%%%%%%%%%%%%%%%%%%%%%%%%%%%
\begin{proof} 
The function $h$ defined by  $h(x):= g(x)+g(-x)$ for $x\in[-a,\,a]$ is
even and convex. Hence on $[0,a]$, $h$ is increasing and $f$ is
decreasing, so that the Chebyshev inequality
obtained by integrating  $(f(x)-f(y))(h(x) -h(y))\leq 0$ over
$[0,a]\times[0,a]$, see \citet[Chapter IX]{MPF93} for references,
yields  $ \int_{-a}^a fg= \int_0^a fh  \leq \frac
1a \int_0^a f \int_0^a h = \frac 1{2a}\int_{-a}^af\int_{-a}^a g $.
\end{proof}
%%%%%%%%%%%%%%%%%%%%%%%%%%%%%%%%%%%%%%%%%%%%%%%%%%%%%%%%%%%%%%%%%%%%%%
\begin{lemma}\label{Wallis ineq}\slshape
For $\lambda\in{]0,\infty[}$, we have 
$\int_0^{\pi/2} \sin^\lambda(t) \,\dd t =
\int_0^{\pi/2} \cos^\lambda(t) \,\dd t< \sqrt{\pi/(2\lambda)}$.
\end{lemma}
%%%%%%%%%%%%%%%%%%%%%%%%%%%%%%%%%%%%%%%%%%%%%%%%%%%%%%%%%%%%%%%%%%%%%%
\begin{proof} 
For $t\in{]0,\pi/2[}$, we have  $\cos(t) =\exp\big(-\int_0^t\tan(u)\,
\dd u\big) < \exp(-t^2/2)$, since $\tan(u)>u$,  so that the
second integral in the claim is  $<$ $\int_0^\infty \exp(-\lambda
t^2/2)\,\dd t$.
\end{proof}
%%%%%%%%%%%%%%%%%%%%%%%%%%%%%%%%%%%%%%%%%%%%%%%%%%%%%%%%%%%%%%%%%%%%%%
\begin{proof}[Proof of the theorem.]
Since the characteristic function $\widehat{\Uniform}_\ell$ of 
$\Uniform_\ell$ is given by 
\[ 
 \widehat{\Uniform}_\ell(t) 
  \,=\,  \frac 1\ell\sum_{k=0}^{\ell-1}\ee^{\ii kt}
  \,=\,   \frac {\ee^{\ii\ell t}-1}{\ell\,
  (\ee^{\ii t} -1)}
  \,=\, \frac {\sin(\ell t/2)}{\ell\,\sin(t/2)}
  \ee^{\ii(\ell-1)t/2}  \quad\qquad (t\in \R)
\]
we get by Fourier inversion for $k\in \Z$
\begin{eqnarray*}
\uniform_\ell^{*n}(k)
&=&\frac{1}{2\pii}
   \int_{-\pii}^\pii 
   \big(\widehat{\Uniform}_\ell(t)
   \big)^n\ee^{-\ii kt}\,\dd t\\
&=&\frac{1}{2\pii}\int_{-\pii}^\pii 
   \Big(\frac{\sin(\ell t/2)}{\ell \sin(t/2)}\Big)^n\exp\Big(\ii
   \Big(\frac{n(\ell -1)}{2}-k\Big) t\Big)\,\dd t\\
&=&\frac{2}{\pii}\int_{0}^{\pii/2} \Big(\frac{\sin(\ell t)
   }{\ell \sin t}\Big)^n\cos((n(\ell -1)-2k)t)\,\dd t
\end{eqnarray*}
Using equality  \eqref{WhereUlnMax}, we get   
\[
c_{\ell,n} \,=\, \frac{2}{\pii}\int_0^{\pii/2} 
\Big(\frac{\sin(\ell t)}{\ell\sin t}\Big)^n\,\cos(\alpha t)\,
\dd t
\,=\, \frac{2}{\pii}\int_0^{\pii/\ell}\,+\,\frac{2}{
\pii}\int_{\pii/\ell}^{\pii/2} 
\,=:\, I_1\,+\,I_2
\]
with 
\[
\alpha\,:=\,n(\ell-1)-2\Big\lfloor \frac{n(\ell-1)}{2}\Big\rfloor
\,\in\, \{0,1\}
\]

To bound $I_1$, we recall the  power series expansion  $x/\tan(x) = 1-
\sum_{k=1}^\infty a_k x^{2k}$  for $|x|<\pii$ with $a_k >0$ for
$k\in \N$,  $a_1=1/3$, and $a_2 = 1/45$, see 
e.g.~\citet[pp.~75--77]{bur79}. With  $b_k := a_k/(2k)$ we get by a 
termwise integration
\[
 -\log  \Big(\frac{\sin x}{x}\Big)
 \,=\, \int_0^x \big(\frac 1y - \frac 1{\tan(y)} \big)\,
 \dd y
 \,=\,\sum_{k=1}^\infty b_k x^{2k}
 \quad\qquad (|x| <\pi)
\]
with $b_k >0$ for $k\in \N$, $b_1= 1/6$, and $b_2 =  1/{180}$. Hence,
for $t\in{]0,\pii/\ell[}$ and with $x:=\sqrt{(\ell^2-1)n/3}\,t$, we
have 
\begin{eqnarray*}
\Big(\frac{\sin(\ell t)}{\ell\sin t}\Big)^n
&=&\exp\Big(n\,\Big(\log\big(\frac{\sin(\ell t)}{\ell t}\big)
-\log\big(\frac{\sin t}{t}\big)\Big)\Big)\\
&=& \exp\Big(-n\sum_{k=1}^\infty b_k
(\ell^{2k}-1)\,t^{2k}\Big) \\
&\leq&\exp\Big(-\frac{n}{6}(\ell^2-1)t^2-\frac{n}{180}
(\ell^4-1)t^4\Big)\\
&\leq&\ee^{-x^2/2}\exp\Big(-\frac{x^4}{20n}\Big)
\qquad\qquad[\text{by $\ell^4-1\ge (\ell^2-1)^2$}]
\end{eqnarray*}
so that, using also $\cos(\alpha t) \leq 1$
and  $\ee^{-y}\leq 1-y+y^2/2$  for $y\in{[0,\,\infty[}$,
\begin{eqnarray}   \label{IneqI1}
\sqrt{\frac{\pii (\ell^2-1)n}{6}}\,I_1&\leq& 
\sqrt{\frac{2(\ell^2-1)n}{3\pii}}\int_{0}^{\pii/\ell}\Big(
\frac{\sin(\ell t)}{\ell\sin t}\Big)^n\,\dd t\\
&\leq&\int_0^{\pii \ell^{-1}
\sqrt{(\ell^2-1)n/3}}\frac{2\ee^{-x^2/2}}{\sqrt{2\pii}}
\exp\Big(-\frac{x^4}{20\,n}\Big)\,\dd x\nonumber\\
&\leq&\int_0^\infty
\frac{2\ee^{-x^2/2}}{\sqrt{2\pii}}
\Big(1-\frac{x^4}{20\,n}+\frac{x^8}{800\,n^2}\Big)\,\dd x\nonumber\\
&=&1-\frac{3}{20\,n}+\frac{21}{160\,n^2}\nonumber
\end{eqnarray}

Now let us bound 
\[
I_2 \,=\, \frac{2}{\pii}\int_{\pii/\ell}^{\pii/2} 
\Big(\frac{\sin(\ell t)}{\ell\sin t}\Big)^n\,\cos(\alpha t)\,\dd t 
\,=\,\int_{\pii}^{\ell\pii/2}\sin^n(t)\,h(t)\,\dd t
\]
where
\[
h(t) \,:=\, \frac {2\cos(\alpha t/\ell)}{\pii\ell\,
 \big(\ell\sin(t/\ell)\big)^n }
 \qquad\qquad (t\in{]0,\ell\pii[}\,)
\]
If $n$ is odd, then 
with $m:= \ell/2$ if $\ell$ is even, $m:=(\ell-1)/2$ if
$\ell \equiv 3 \pmod 4$,  and $m := (\ell +1)/2$ if 
$\ell \equiv 5 \pmod 4$, we get
\begin{equation} \label{IenqI2nodd}
 I_2 \,\leq\, \int_\pii^{m\pii} \sin^n(t)\,h(t)\,\dd t
 \,=\, \int_0^\pii\sin^n(t)\sum_{j=1}^{m-1}(-1)^j h(t+j\pii) 
 \,\dd t  \,\leq\, 0
\end{equation}
since $h$ is positive and decreasing. 
If $n$ is even, then we use $\cos x \leq 1$ and  
$\sin x \ge 2x/\pii$ for $x\in[0,\pi/2]$ 
to get $h(t) \leq \frac 1\ell  (\pii/2)^{n-1}/t^n$ and hence
\begin{eqnarray}\label{IneqI2neven}
I_2
 &\leq & \frac{1}{\ell}\Big(\frac{\pii}{2}\Big)^{n-1}
\sum_{k=1}^\infty
\int_{k\pii}^{(k+1)\pii}\sin^n( t)\,\frac{1}{t^n}\,\dd t\\
  &\leq&  \frac{1}{\ell}\Big(\frac{\pii}{2}\Big)^{n-1}
\sum_{k=1}^\infty
\frac 1\pi\int\limits_{k\pii}^{(k+1)\pii} \sin^n (t)\,\dd t
\int\limits_{k\pii}^{(k+1)\pii}\frac{ \dd t  }{t^n}
\nonumber \\  
&& \qquad\qquad \qquad \qquad  \qquad \quad  
 [\text{Lemma~\ref{l3}, $t=x +  (k+\frac 12)\pi $}] \nonumber \\  
&=&\frac{1}{\pii\,\ell\,(n-1)\,2^{n-1}}
\int_{0}^{\pii} \sin^n (t)\,\dd t \nonumber \\
&\leq&  \frac{1}{\ell\,(n-1)\,2^{n-1}}\sqrt{\frac{2}{
\pii n}}  \quad \qquad\qquad\qquad\qquad 
 [\text{by Lemma \ref{Wallis ineq}}]   \nonumber
\end{eqnarray}

Combining our estimates from
\eqref{IneqI1}, \eqref{IenqI2nodd}, \eqref{IneqI2neven} 
and using
$\sqrt{\ell^2-1}< \ell$, we obtain 
\begin{equation}\label{ee88}
\sqrt{\frac{\pii (\ell^2-1)n}{6}}\,c_{\ell,n} \,\leq\, 
1-\frac{3}{20\,n}+\frac{21}{160\,n^2}+
\frac{\1^{}_{2\N}(n)}{\sqrt{3}\,(n-1)\,2^{n-1}}
\,=:\, d_n
\end{equation}
for all $\ell, n\in \N $ with $\ell \ge 2$. 
For n odd, we use $n^2 \ge n$ to get
$d_n -1 \leq \frac 1{n}(-\frac 3{20}+\frac {21}{160})<0$.
For $n$ even with  $n \neq 2$, we use $n^2 \ge 4n$
and $(n-1)\,2^{n-1}\ge 6n  $ in \eqref{ee88} to get 
\[
 d_n -1 \,\leq\, \frac {1}{n}\, \Big(- \frac{3}{20} +
\frac{21}{160}\cdot \frac 14  + \frac{1}{6\sqrt{3}}\Big) 
\,=\,\frac{1}{2n}\,\Big(\frac{1}{3\sqrt{3}}-\frac{15}{64}\Big) 
\,<\, 0
\]
Thus  for $n\neq 2$, we  have $d_n <1 $,  and hence inequality 
\eqref{MainInequality}. For   $n=2$, the claim  of the theorem follows
from \eqref{cln for n=1,2}.
\end{proof}
%%%%%%%%%%%%%%%%%%%%%%%%%%%%%%%%%%%%%%%%%%%%%%%%%%%%%%%%%%%%%%%%%%%%%%
\section*{Acknowledgement}
We thank Jannis Dimitriadis for thoroughly reading an earlier draft 
of this note.
%%%%%%%%%%%%%%%%%%%%%%%%%%%%%%%%%%%%%%%%%%%%%%%%%%%%%%%%%%%%%%%%%%%%%%

%%%%%%%%%%%%%%%%%%%%%%%%%%%%%%%%%%%%%%%%%%%%%%%%%%%%%%%%%%%%%%%%%%%%%%
%%%%%%%%%%%%%%%%%%%%%%%%%%%%%%%%%%%%%%%%%%%%%%%%%%%%%%%%%%%%%%%%%%%%%%
\end{document}